\documentclass[11pt,reqno]{amsart}
\usepackage{color}

\newtheorem{theorem}{Theorem}
\newtheorem{proposition}
{Proposition}
\newtheorem{corollary}[proposition]{Corollary}
\newtheorem{remark}[proposition]{Remark}
\newtheorem{lemma}[proposition]{Lemma}
\newfont{\bb}{msbm10 at 12pt}

\def\pf{{\textit {Proof :} }}

\def\End{{\rm End}}    \def\Id{{\rm Id}}
\def\C{{\mathbb C}}

\def\Sps{\mathcal{S}}

\def\bS{{\mathbb S}}

\def\E{\mathcal{E}}
\def\So{{\mathbb{S} \Omega}}
\def\<{\langle}     
\def\>{\rangle}

\def\UB{\mathbb{B}^{n+1}}
\def\Spin{{\rm Spin}}
\def\SO{{\rm SO}}

\newcommand{\bal}{\begin{align}}      \newcommand{\eal}{\end{align}}
\newcommand{\ba}{\begin{array}}      \newcommand{\ea}{\end{array}}
\newcommand{\bc}{\begin{center}}     \newcommand{\ec}{\end{center}}
\newcommand{\be}{\begin{enumerate}}  \newcommand{\ee}{\end{enumerate}}
\newcommand{\beQ}{\begin{eqnarray*}} \newcommand{\eeQ}{\end{eqnarray*}}
\newcommand{\bi}{\begin{itemize}}    \newcommand{\ei}{\end{itemize}}
\newcommand{\bt}{\begin{tabular}}    \newcommand{\et}{\end{tabular}}
\newcommand{\bdm}{\begin{displaymath}} \newcommand{\edm}{\end{displaymath}}

    \newcommand{\sm}{\bS\!\!\!/\,\!}
\newcommand{\D}{D\!\!\!\!/\,}
\newcommand{\nb}{\nabla\!\!\!\!/\,}
\newcommand{\mult}{\gamma\!\!\!/}
\newcommand{\Eb}{\E\!\!\!/}

\newcommand{\lto}{\ensuremath{\longrightarrow}}

\def\qed{\hfill{q.e.d.}\smallskip\smallskip}

\begin{document}

\title[On a quasi-local mass]{A holographic principle for the existence of imaginary Killing spinors}

\author{Oussama Hijazi}
\address[Oussama Hijazi]{Institut {\'E}lie Cartan de Lorraine,
Universit{\'e} de Lorraine, Nancy,
B.P. 239,
54506 Vand\oe uvre-L{\`e}s-Nancy Cedex, France.}
\email{Oussama.Hijazi@univ-lorraine.fr}

\author{Sebasti{\'a}n Montiel}
\address[Sebasti{\'a}n  Montiel]{Departamento de Geometr{\'\i}a y Topolog{\'\i}a,
Universidad de Granada,
18071 Granada,  Spain.}
\email{smontiel@ugr.es}

\author{Simon Raulot}
\address[Simon Raulot]{Laboratoire de Math\'ematiques R. Salem
UMR $6085$ CNRS-Universit\'e de Rouen
Avenue de l'Universit\'e, BP.$12$
Technop\^ole du Madrillet
$76801$ Saint-\'Etienne-du-Rouvray, France.}
\email{simon.raulot@univ-rouen.fr}

\begin{abstract}
Suppose that $\Sigma=\partial\Omega$ is the $n$-dimensional boundary, with positive (inward) mean curvature $H$, of a connected  compact $(n+1)$-dimensional Riemannian spin manifold $(\Omega^{n+1},g)$ whose scalar curvature $R\ge -n(n+1)k^2$,  for some $k>0$. If $\Sigma$ admits an isometric and isospin immersion $F$ into the hyperbolic space ${\mathbb{H}^{n+1}_{-k^2}}$,
 we define a quasi-local mass and prove its positivity as well as the associated rigidity statement. The proof is based on a holographic principle for the existence of an imaginary Killing spinor. For $n=2$, we also show that its limit, for coordinate spheres in an Asymptotically Hyperbolic (AH) manifold, is the mass of the (AH) manifold. 	 
\end{abstract}



\thanks{The second author was partially  
supported by a Spanish MEC-FEDER grant No. MTM2007-61775}

\date{\today}   
\maketitle 
\pagenumbering{arabic}
\begin{center}
 \address{\em  D\'edi\'e \`a {\rm Paul Gauduchon} en t\'emoignage de notre
reconnaissance et amiti\'e.}
\end{center}

\section{Introduction}



The Positive Mass Theorem (PMT) states that for a complete
asymptotically flat manifold which, near each end, behaves like the Euclidean
space at infinity and whose scalar curvature is nonnegative,
then its ADM mass of each end is non-negative. Moreover,
if the ADM mass of one end is zero, then the manifold is the
Euclidean space. The PMT was proved by Schoen
and Yau \cite{SY1, SY2} using minimal surface techniques. Later on, Witten
\cite{Wi} gave an elegant and simple  proof of the PMT for spin manifolds.
Since then, spinors has been successfully used
 to prove Positive Mass type theorems (see for example \cite{AD, He, CH, LY1, LY2, Wa1, ST1, HM1}). \\

In this spirit,  Wang and
Yau \cite{WY} introduced a quasi-local mass for $3$-dimensional manifolds
with boundary whose scalar curvature is bounded from below by a
negative constant. Again, using spinorial methods, they  proved that
this mass is non-negative. Shi and Tam \cite{ST2}  proved a
similar result  but with a simpler and more explicit definition of the mass. More precisely:

 
\begin{theorem}\label{ST} \cite{ST2}
Let $(\Omega^3,g)$ be a compact $3$-dimensional orientable Riemannian manifold with smooth boundary $\Sigma$. Assume that :
\begin{enumerate}
\item The scalar curvature $R$ of $(\Omega,g)$ satisfies $R\geq -6k^2$ for some $k>0$;
\item The boundary $\Sigma$ is a topological sphere with Gauss curvature \hfill\break $K>-k^2$ and  mean curvature $H>0$ (so that
$\Sigma$ can be isometrically embedded into $\mathbb{H}^3_{-k^2}$, the Hyperbolic space of 
 constant curvature $-k^2$, with mean curvature $H_0$). 
\end{enumerate}
Then, the energy-momentum vector
\begin{eqnarray*}
\mathcal{M}_\alpha:=\int_\Sigma(H_0-H){\mathbf{W}_\alpha}\,d\Sigma\in\mathbb{R}^{3,1}
\end{eqnarray*}
is future directed non-spacelike or zero, where ${\mathbf{W}_\alpha}=(x_1,x_2,x_3,\alpha t)$ with 
\begin{eqnarray}\label{alpha}
\alpha=\coth R_1+\frac{1}{\sinh R_1}\Big(\frac{\sinh^2R_2}{\sinh^2 R_1}-1\Big)^{\frac{1}{2}}>1
\end{eqnarray}
an explicit  constant depending on the intrinsic geometry of $\Sigma$ and ${\bf X} :={\mathbf{W}_1}=(x_1,x_2,x_3,t)$ is the position vector in $\mathbb{R}^{3,1}$. Moreover, if there exists a future directed null vector $\zeta\in\mathbb{R}^{3,1}$ such that:
\begin{eqnarray*}
\<\mathcal{M}_\alpha,\zeta\>_{\mathbb{R}^{3,1}}=0,
\end{eqnarray*}
then $(\Omega^3,g)$ is a domain in $\mathbb{H}^3_{-k^2}$.
\end{theorem}
The statement of this result needs some explanation. First, from \cite{P} and \cite{DCW}, as mentioned, the assumptions on the boundary $\Sigma$ ensure the existence of an isometric embedding of $\Sigma$ into the hyperbolic space $\mathbb{H}^3_{-k^2}$ as a convex surface which bounds a domain $D$ in $\mathbb{H}^3_{-k^2}$. Moreover, this embedding is unique up to an isometry of $\mathbb{H}^3_{-k^2}$. Here $H_0$ denotes the mean curvature of this embedding and $R_1$ and $R_2$ are two positive real numbers such that $B_o(R_1)\subset D \subset B_o(R_2)$ in $\mathbb{H}^3_{-k^2}$ where $B_o(r)$ is the geodesic ball of radius $r>0$ and center $o=(0,0,0,1/k)$. This result has been recently generalized by Kwong \cite{K}. Namely, he proves:
\begin{theorem}\label{K} \cite{K}
For $n\geq 2$, let $(\Omega^{n+1},g)$ be a compact spin $(n+1)$-dimensional manifold with smooth boundary $\Sigma$. Assume that:
\begin{enumerate}
\item The scalar curvature $R$ of $\Omega$ satisfies $R\geq -n(n+1)k^2$ for some $k>0$,
\item The boundary $\Sigma$ is topologically an $n$-sphere with sectional curvature $K>-k^2$, mean curvature $H>0$ and that $\Sigma$ can be isometrically embedded uniquely into ${\mathbb{H}^{n+1}_{-k^2}}$ with mean curvature $H_0$. 
\end{enumerate}
Then, there is a future time-like vector-valued function ${\bf W}_\alpha$ on $\Sigma$ such that the energy-momentum vector:
\begin{eqnarray*}
\mathcal{M}_\alpha:=\int_\Sigma(H_0-H){\bf W}_\alpha\,d\Sigma\in\mathbb{R}^{n+1,1}
\end{eqnarray*}
is future non-spacelike. Here ${\bf W}_\alpha =(x_1,x_2,\cdots,x_{n+1},\alpha t)$ for some $\alpha>1$ and ${\bf X} := {\bf W}_1(x_1,x_2,\cdots,x_{n+1},t)\in{\mathbb{H}^{n+1}_{-k^2}}\subset\mathbb{R}^{n+1,1}$ is the position vector of the embedding of $\Sigma$.
\end{theorem}
In the general case,  the constant $\alpha$ is still explicitly given by (\ref{alpha}). It is conjectured, and verified for $n=2$ in certain cases (see \cite{ST2}), that Theorems \ref{ST} and \ref{K} should hold for $\alpha=1$. A key ingredient in the proof of these two results is a generalization of the Positive Mass Theorem for (AH) manifolds (see Section \ref{AHset}).

In this paper, we make use of another approach, developed in \cite{HM1}, to establish a holographic principle\footnote{By holographic principle we mean the property which states that the description of a manifold with boundary 
can be thought of as encoded on the boundary.} for the existence of imaginary Killing spinors on Dirac bundles (See Section \ref{Sp-Ch}) in order to generalize the above results in several directions. 
Namely, we modify the curvature term in the definition of $\mathcal{M}_\alpha$ to precisely define a energy-momentum vector field $\mathbf{E}(\Sigma)$ in terms of 
${\bf X}$. In particular, our expression depends only on the metric and the embedding of $\Sigma$ and is thus independent of the particular manifold $\Omega$. It could be considered as a possible new definition of a quasi-local mass since it has the desirable {\it non negativity} and {\it rigidity} properties as shown in Theorems \ref{general} and \ref{3D}. Moreover, these statements hold in a more general setup.

In fact, we have:
\begin{theorem}\label{general}
Let $(\Omega^{n+1},g)$ be a compact, connected $(n+1)$-dimensional Riemannian spin manifold with smooth boundary $\Sigma$. Assume that
\begin{enumerate}
\item
The scalar curvature $R$ of $\Omega$ satisfies $R\geq -n(n+1)k^2$ for some $k>0$;
\item
The boundary $\Sigma = \partial \Omega$ has  mean curvature $H >0$  and that 
there exists an isometric and isospin immersion $F$ of $\Sigma$ into the hyperbolic space ${\mathbb{H}}^{n+1}_{-k^2}$ with mean curvature $H_0$.
\end{enumerate}
 Then, the energy-momentum vector defined by 
\begin{eqnarray}\label{E-M-V}
\mathbf{E}(\Sigma):=\int_\Sigma\Big(\frac{H_0^2-H^2}{H}\Big){\bf X}\,d\Sigma\in\mathbb{R}^{n+1,1}
\end{eqnarray}
is timelike future directed or zero (see Theorem \ref{K} for the definition of $X$). Moreover, ${\bf E}(\Sigma)=0$ if and only if $(\Omega^{n+1},g)$ is a domain in $\mathbb{H}^{n+1}_{-k^2}$, $\Sigma$ is connected and the embedding of $\Sigma$ in $\Omega$ and its immersion $F$ in ${\mathbb{H}}^{n+1}_{-k^2}$ are congruent.
\end{theorem}

For $n=2$, since $\Omega$ is automatically spin, we deduce the following:
\begin{theorem}\label{3D}
Let $(\Omega^3,g)$ be a compact, connected $3$-dimensional oriented Riemannian manifold with smooth boundary $\Sigma$. Assume that:
\begin{enumerate}
\item The scalar curvature $R$ of $(\Omega^3,g)$ satisfies $R\geq -6k^2$ for some $k>0$;
\item The boundary $\Sigma$ is a topological sphere with Gauss curvature \hfill\break $K>-k^2$ and with  mean curvature $H >0$.
\end{enumerate}
Then, the energy-momentum vector given by 
$\mathbf{E}(\Sigma)
\in\mathbb{R}^{3,1}$
is timelike future directed or zero. Moreover, $\mathbf{E}(\Sigma)=0$ if and only if $(\Omega^3,g)$ is a domain in $\mathbb{H}^{3}_{-k^2}$ and $\Sigma$ is connected.
\end{theorem}

Note that in the rigidity part of this result, the embedding of $\Sigma$ in $\Omega$ and its immersion in ${\mathbb{H}}^{3}_{-k^2}$ are automatically congruent because of the uniqueness of the embedding of $\Sigma$ in ${\mathbb{H}}^{3}_{-k^2}$.\\

For simplicity, we will only prove the case $k=1$. The general case is obtained by a homothetic change of the metric.


\section{Geometric and Analytic preliminaries}


The aim of this section is to introduce the general geometrical spinorial setting and the basic  analytical tools needed to establish the above mentioned results.


\subsection{The geometric setting}\label{G-S}


In the following, we consider a compact and connected Riemannian spin $(n+1)$-dimensional manifold $(\Omega^{n+1},g)$ with smooth boundary $\Sigma:=\partial\Omega$. The Riemannian structure on $\Omega$ induces a Riemannian metric on $\Sigma$, also denoted by $g$, whose Levi-Civita connection $\nabla^\Sigma$ satisfies the Riemannian Gauss formula
\begin{eqnarray}\label{f-g-r}
\nabla^\Sigma_XY=\nabla^\Omega_XY-g\big(A(X),Y\big)N
\end{eqnarray}
for all $X,Y\in\Gamma(T\Sigma)$. Here $\nabla^\Omega$ is the Levi-Civita connection on $\Omega$, $N$ the unit inner normal vector field to $\Sigma$ and $A$ is the Weingarten map defined by $A(X)=-\nabla^\Omega_XN$, for $X\in\Gamma(T\Sigma)$. Since $\Omega$ is spin, there exists a pair $(\Spin(\Omega),\eta)$ where $\Spin(\Omega)$ is a $\Spin_{n+1}$-principal fiber bundle over $\Omega$ and $\eta$ is a $2$-fold covering of the $\SO_{n+1}$-principal bundle $\SO(\Omega)$ of $g$-orthonormal frames such that 
\begin{eqnarray*}
\forall u\in\Spin(\Omega),\,\,\forall a\in\Spin_{n+1},\qquad \eta(ua)=\eta(u)\rho(a)
\end{eqnarray*}
where $ua$ denotes the right action of $\Spin_{n+1}$ on $\Spin(\Omega)$ and $\rho$ is the $2$-fold covering of the special orthogonal group $\SO_{n+1}$ by $\Spin_{n+1}$. Note that since $\Omega$ is oriented, it induces an orientation on the boundary, hence $\Sigma$ is automatically spin. Indeed, via the inclusion map $\SO(\Sigma)\hookrightarrow\SO(\Omega)$, we can define the pulled-back bundle $\Spin(\Omega)_{|\Sigma}$, which gives a spin structure on $\Sigma$ denoted by $\Spin(\Sigma)$. Recall that on $\Omega$, we define the spinor bundle $\So$, a rank $2^{\left[\frac{n+1}{2}\right]}$ complex vector bundle, by 
$$\So:=\Spin(\Omega)\times_{\gamma_{n+1}}\mathcal{S}_{n+1}$$
where $\gamma_{n+1}$ is the restriction to $\Spin_{n+1}$ of an irreducible complex representation of the complex Clifford algebra $\C l_{n+1}$. This representation provides a left Clifford module 
\begin{eqnarray}
\gamma^\Omega:\C l(\Omega)\longrightarrow\End_\C(\So)
\end{eqnarray}
which is a fiber preserving algebra morphism. Then $\So$ becomes a bundle of complex left modules over the Clifford bundle $\C l(\Omega)$. In particular, $\So$ is a complex {\it Dirac bundle} in the sense of \cite{LM}, i.e.,  there exists on $\mathbb{S}\Omega$:
\begin{itemize}
\item[-] a Hermitian scalar product $\<\,,\,\>_{\Omega}$,
\item[-] a spin Levi-Civita connection $\nabla^\Omega$ acting on sections of $\So$
\end{itemize}
such that
\begin{itemize}
\item[-] the Clifford multiplication by tangent vector fields is skew-Hermitian: 
\begin{eqnarray}\label{s-h}
\<\gamma^\Omega(X)\psi , \varphi\>=-\<\psi , \gamma^\Omega(X)\varphi\>,\end{eqnarray}
\item[-] the covariant derivative $\nabla^\Omega$ is a module derivation, that is 
\begin{eqnarray}\label{m-d}
\nabla^\Omega_X\big(\gamma^{\Omega}(Y)\psi \big)=\gamma^{\Omega}(\nabla^\Omega_XY)\psi +\gamma^\Omega(Y)\nabla^\Omega_X\psi ,
\end{eqnarray}
\item[-] the covariant derivative $\nabla^\Omega$ is compatible with the Hermitian scalar product, that is 
\begin{eqnarray}\label{c-d}
X\<\psi, \varphi\>=\<\nabla^\Omega_X\psi, \varphi\>+\< \psi ,\nabla^\Omega_X \varphi\>
\end{eqnarray}
\end{itemize}
for all $X,Y\in\Gamma(T\Omega)$ and $\psi, \varphi\in\Gamma(\So)$. The Dirac operator $D^\Omega$ on $\So$ is the first order 
elliptic differential operator locally given by 
$$
D^\Omega=\sum_{i=1}^{n+1}\gamma^\Omega(e_i)\nabla^\Omega_{e_i},
$$
where $\{e_1,\dots,e_{n+1}\}$ is a local orthonormal frame of $T\Omega$. As mentioned above, the boundary is naturally endowed with a spin structure and the group $\Spin_n\subset\C l_n^0$ (the even part of the Clifford algebra) acts on the restricted bundle $\Spin(\Omega)_{|\Sigma}$ via the map $\iota$ defined by 
$$
\begin{array}{r@{}ccl}
\iota \colon  & \C l_n & \lto & \C l_{n+1}^0\subset\C l_{n+1} \\[.05cm]  & e_j &\longmapsto  & e_j\cdot N.
\end{array}
$$
where the dot is the multiplicative structure of the Clifford algebra. Hence we have that the restriction 
\begin{eqnarray*}
\sm\Sigma:=\So_{|\Sigma}=\Spin(\Sigma)\times_{\gamma_{n+1}\circ\iota}\mathcal{S}_{n+1}
\end{eqnarray*}
is a left module over $\C l(\Sigma)$ with Clifford multiplication
\begin{eqnarray*}
\mult^\Sigma:\C l(\Sigma)\longrightarrow\End_\C(\sm\Sigma)
\end{eqnarray*}
given by $\mult^\Sigma=\gamma^\Omega\circ\iota$, that is 
\begin{eqnarray}\label{e-c-m}
\mult^\Sigma(X)\psi=\gamma^\Omega(X)\gamma^\Omega(N)\psi
\end{eqnarray}
for every $\psi\in\Gamma(\sm\Sigma)$ and $X\in\Gamma(T\Sigma)$. Consider on $\sm\Sigma$ the Hermitian
metric $\<\,,\,\>_{\Omega}$ induced from that of $\So$. This metric immediately satisfies
the compatibility condition (\ref{s-h}) if one puts on $\Sigma$ the Riemannian metric induced from $\Omega$ and the extrinsic Clifford multiplication $\mult^\Sigma$ defined in (\ref{e-c-m}). Now the Gauss formula (\ref{f-g-r}) implies that the spin connection $\nb^\Sigma$ on $\sm\Sigma$ is given by the following spinorial Gauss formula
\begin{eqnarray}\label{s-g-f}
\nb^\Sigma_X\psi=\nabla^\Omega_X\psi-\frac{1}{2}\mult^\Sigma(AX)\psi
\end{eqnarray}
for every $\psi\in\Gamma(\sm\Sigma)$ and $X\in\Gamma(T\Sigma)$. The extrinsic Dirac operator\hfill\break $\D^\Sigma:=\mult^\Sigma\circ\nb^\Sigma$ on $\Sigma$ defines a first order elliptic operator acting on sections of $\sm\Sigma$. By (\ref{s-g-f}), for any spinor field $\psi\in\Gamma(\So)$, we have 
\begin{equation} \label{D-D}
{\D}^\Sigma\psi=\sum_{j=1}^n\mult^\Sigma(e_j)\nb^\Sigma_{e_j}\psi
=\frac{n}{2}H\psi-\gamma^\Omega(N)D^\Omega\psi-\nabla^\Omega_N\psi,
\end{equation}
and
\begin{equation}\label{D-commutes}
\D^\Sigma\big(\gamma^\Omega(N)\psi\big)=-\gamma^\Omega(N)\D^\Sigma\psi
\end{equation}
where $\{e_1,\dots,e_n\}$ is a local orthonormal frame of $T\Sigma$ and $H=\frac{1}{n}\hbox{trace\,}A$ is the mean curvature of $\Sigma$ in $\Omega$. On the other hand, $\Sigma$ has also an {\it intrinsic} spinor bundle defined from its spin structure and an irreducible representation of $\C l_n$. More precisely, the complex vector bundle of rank $2^{[\frac{n}{2}]}$, defined by
\begin{eqnarray*}
\bS\Sigma:=\Spin(\Sigma)\times_{\gamma_n}\Sps_n .
\end{eqnarray*}
This is also a Dirac bundle over $\Sigma$ with a Clifford multiplication $\gamma^\Sigma$, a spin Levi-Civita connection $\nabla^\Sigma$ and a  Hermitian scalar product satisfying the properties (\ref{s-h}), (\ref{m-d}) and (\ref{c-d}) on $\Sigma$. Moreover, the intrinsic Dirac operator on $\Sigma$ is then defined by $D^\Sigma:=\gamma^\Sigma\circ\nabla^\Sigma$. As we shall see in the next section, there are natural identifications between intrinsic and extrinsic spinor bundles over $\Sigma$ (see \cite{Bu, Tr, Baer, HMZ3, HMR2}  for more details).

\subsection{Dirac bundles and chirality operator}\label{Sp-Ch}

The important fact now is to consider bundles on which a chirality operator is defined. Recall that a chirality operator $\omega$ on a Dirac bundle $(\E\Omega,\gamma,\nabla,\<\,,\,\>)$ is an endomorphism
\begin{eqnarray*}
\omega:\Gamma(\E\Omega)\longrightarrow\Gamma(\E\Omega)
\end{eqnarray*}
such that
\begin{eqnarray}
\label{c-o-1}\omega^2=\Id_{\E\Omega}, \qquad& \< \omega\Psi , \omega\Phi\>=\<\Psi, \Phi\>,\\
\label{c-o-2}\omega(\gamma(X)\Psi )=-\gamma(X)\omega\Psi , \qquad& \nabla_X(\omega\Psi )=\omega(\nabla_X\Psi ),
\end{eqnarray}
for all $X\in\Gamma(T\Omega)$ and $\Psi ,\Phi\in\Gamma(\E\Omega)$.
\vspace{0.2cm}
In the following, we consider the vector bundle given by
$$
\E\Omega:=\left\lbrace
\begin{array}{ll}
\So\quad & {\rm if}\,\, n+1=2m,\\
\So\oplus\So\quad  & {\rm if }\,\, n+1=2m+1,
\end{array}
\right.
$$
on which a Clifford multiplication $\gamma$ and a linear connection $\nabla$ are defined by
\begin{equation}\label{c-m-e}
\gamma=\left\lbrace
\begin{array}{ll}
\gamma^\Omega\quad & {\rm if }\,\, n+1=2m\\ \\
\gamma^\Omega\oplus-\gamma^\Omega=
\begin{pmatrix}
\gamma^\Omega & 0\\
0 & -\gamma^\Omega
\end{pmatrix}\quad  & {\rm if }\,\, n+1=2m+1 
\end{array}
\right.
\end{equation}
and
\begin{equation}\label{l-c-e}
\nabla=\left\lbrace
\begin{array}{ll}
\nabla^\Omega\quad & {\rm if }\,\, n+1=2m\\ \\
\nabla^\Omega\oplus\nabla^\Omega=
\begin{pmatrix}
\nabla^\Omega\quad & 0\\
0 & \nabla^\Omega
\end{pmatrix}  
& {\rm if }\,\, n+1=2m+1.
\end{array}
\right.
\end{equation}
Finally, $\<\,,\,\>$ denotes the Hermitian scalar product given by $\<\,,\,\>_{\Omega}$ for $n$ odd and by
\begin{eqnarray}\label{h-s-o}
\<\Psi\ , \Phi\>:=\<\psi_1 , \varphi_1\>_\Omega+\<\psi_2 , \varphi_2\>_\Omega
\end{eqnarray}
for $n$ even, for any $\Psi=(\psi_1,\psi_2)$, $\Phi=(\varphi_1,\varphi_2)\in\Gamma(\E\Omega)$. The Dirac-type operator acting on sections of $\E\Omega$ and defined by $D:=\gamma\circ\nabla$ is explicitly given by
$$
D=\left\lbrace
\begin{array}{ll}
D^\Omega\quad & {\rm if }\,\, n+1=2m\\ \\
D^\Omega\oplus-D^\Omega=
\begin{pmatrix}
D^\Omega & 0\\
0 & -D^\Omega
\end{pmatrix}\quad & {\rm if }\,\, n+1=2m+1.
\end{array}
\right.
$$
Let us examine,  in more details, this bundle and its restriction to $\Sigma$:

\vspace{0.2cm}

\noindent{\it The even dimensional case}

If $n+1=2m$, the vector bundle $\E\Omega$ is the spinor bundle $\So$. In this situation, it is well-known that the Clifford multiplication $\omega:=\gamma(\omega_{n+1}^\C)$ by the complex volume element 
\begin{eqnarray*}
\omega^\C_{n+1}:=i^me_1\cdot...\cdot e_{n+1}
\end{eqnarray*}
defines a chirality operator on $\E\Omega$. Moreover, the spinor bundle splits into
\begin{eqnarray}\label{dec-chiral}
\E\Omega=\So=\mathbb{S}^+\Omega\oplus\mathbb{S}^-\Omega
\end{eqnarray}
where $\mathbb{S}^\pm\Omega$ are the $\pm 1$-eigenspace of the endomorphism $\omega$. On the other hand, from algebraic considerations (see \cite{HMZ1} or \cite{HMR2} for example) the restricted spinor bundle 
\begin{eqnarray*}
\Eb:=\E\Omega_{|\Sigma}=\So_{|\Sigma}=\sm\Sigma
\end{eqnarray*}
can be identified with the  intrinsic data of  $\Sigma$ as follows:
\begin{eqnarray*}
(\sm\Sigma,\mult^{\Sigma},\nb^\Sigma)\cong(\mathbb{S}\Sigma\oplus\mathbb{S}\Sigma,\gamma^\Sigma\oplus-\gamma^\Sigma,\nabla^\Sigma\oplus\nabla^\Sigma). 
\end{eqnarray*}  
In the following, for simplicity we let $(\Eb,\mult,\nb):=(\sm\Sigma,\mult^\Sigma,\nb^\Sigma)$ the extrinsic Dirac bundle over the boundary of the even dimensional Riemannian spin domain $\Omega$. As a consequence of these identifications, we get that the extrinsic Dirac-type operator $\D:=\mult\circ\nb$ of $\Eb$ can be identified with the  extrinsic Dirac operator $\D^\Sigma$ which only depends on the Riemannian and the spin structure of $\Sigma$ since we also have the following identification
\begin{eqnarray*}
\D=D^\Sigma\oplus-D^\Sigma=
\begin{pmatrix}
D^\Sigma & 0\\
0 & -D^\Sigma
\end{pmatrix}.
\end{eqnarray*}
Moreover, a simple but important observation here is that we can also choose the Clifford action of the unit normal $N$ by: 
\begin{eqnarray}\label{a-n-p}
\gamma(N)=-i
\begin{pmatrix}
0 & \Id\\
\Id & 0
\end{pmatrix}
\end{eqnarray}
where the matrix blocks are defined with respect to the chiral decomposition (\ref{dec-chiral}). Then we note that the Dirac-type operator defined for all $\Psi\in\Gamma(\Eb)$ by
\begin{eqnarray*}
\D^\pm\Psi:=\D\Psi\pm\frac{n}{2}i\gamma(N)\Psi=\D^\Sigma\pm\frac{n}{2}i\gamma^\Omega(N)\Psi
\end{eqnarray*}
does not depend on the extrinsic geometry of $\Sigma$ in $\Omega$. Indeed, from the identification of $\D$ and (\ref{a-n-p}), we have 
\begin{eqnarray*}
\D^\pm=
\begin{pmatrix}
D^{\Sigma} & \pm\frac{n}{2}\Id\\
\pm\frac{n}{2}\Id & -D^\Sigma
\end{pmatrix}
\end{eqnarray*}
and it is obvious from this expression that these operators only depend on intrinsic data of $\Sigma$ (more precisely on the spin structure and the induced metric of $\Sigma$).

\vspace{0.2cm}

\noindent{\it The odd dimensional case}

If $n+1=2m+1$, the vector bundle $\E\Omega$ consists of two copies of the  spinor bundle
\begin{eqnarray*}
\E\Omega=\So\oplus\So
\end{eqnarray*}
and its rank on $\C$ is $2^{m+1}$. It is straightforward from the definitions (\ref{c-m-e}), (\ref{l-c-e}) and (\ref{h-s-o}) that the relations (\ref{s-h}), (\ref{m-d}) and (\ref{c-d}) are valid for $\gamma$, $\nabla$ and $\<\,,\,\>$ and thus $(\E\Omega,\gamma,\nabla,\<\,,\,\>)$ defines a Dirac bundle over $\Omega$. In this situation, it is a simple exercise to check that the map 
$$
\begin{array}{r@{}ccl}
\omega \colon  & \Gamma(\E\Omega) &\lto & \Gamma(\E\Omega) \\[.05cm]  
& \Psi=\begin{pmatrix}\psi_1\\ \psi_2\end{pmatrix} &\longmapsto  & \omega\Psi:=\begin{pmatrix}\psi_2\\ \psi_1\end{pmatrix}, 
\end{array}
$$
satisfies the properties (\ref{c-o-1}) and (\ref{c-o-2}) so that it defines a chirality operator on $\E\Omega$. The restriction of $\E\Omega$ to $\Sigma$ is given by
\begin{eqnarray*}
\Eb:=\E\Omega_{|\Sigma}=\sm\Sigma\oplus\sm\Sigma
\end{eqnarray*}
and can be identified with two copies of the intrinsic spinor bundle of $\Sigma$ (see \cite{HMZ1} or \cite{HMR2} for more details). Similarly, the extrinsic spin Levi-Civita connection 
\begin{equation}\label{l-c-oe}
\nb:=\nb^\Sigma\oplus\nb^\Sigma=
\begin{pmatrix}
\nb^\Sigma & 0\\
0 & \nb^\Sigma
\end{pmatrix}
\end{equation} 
as well as its Clifford multiplication 
\begin{equation}\label{c-m-oe}
\mult:=\mult^\Sigma\oplus\mult^\Sigma=
\begin{pmatrix}
\mult^\Sigma & 0\\
0 & \mult^\Sigma
\end{pmatrix}
\end{equation}
are such that the following identifications hold
\begin{eqnarray*}
(\Eb,\mult,\nb)\cong(\sm\Sigma\oplus\sm\Sigma,\gamma^\Sigma\oplus\gamma^\Sigma,\nabla^\Sigma\oplus\nabla^\Sigma).
\end{eqnarray*} 
In particular, these definitions provide a  Dirac bundle structure on $\Eb$. It is also clear from the definitions of $\nabla$, $\nb$ and the spinorial Gauss formula (\ref{s-g-f}) that a similar relation holds between $\nabla$ and $\nb$. The extrinsic Dirac-type operator acting on sections of $\Eb$ is defined as usually by $\D:=\mult\circ\nb$ and by (\ref{l-c-oe}) and (\ref{c-m-oe}), it satisfies: 
\begin{eqnarray*}
\D=
\begin{pmatrix}
D^\Sigma & 0\\
0 & D^\Sigma
\end{pmatrix}.
\end{eqnarray*}
Then we also easily observe that  relations (\ref{D-D}) and (\ref{D-commutes}) hold. Finally, as in the even dimensional case, the operators defined by 
\begin{eqnarray}\label{MEDTO}
\D^\pm:=\D\pm\frac{n}{2}i\gamma(N)
\end{eqnarray}
can be expressed intrinsically with respect to $\Sigma$. Indeed, by (\ref{c-m-e}), we first note that
\begin{eqnarray*}
\D^\pm=
\begin{pmatrix}
D^\Sigma\pm\frac{n}{2}i\gamma^\Omega(N) & 0 \\
0 & D^\Sigma\mp\frac{n}{2}i\gamma^\Omega(N)
\end{pmatrix}.
\end{eqnarray*}
Moreover, since
$$
D^{\Sigma}:\Gamma\big(\mathbb{S}^\pm(\Sigma)\big)\longrightarrow\Gamma\big(\mathbb{S}^\mp(\Sigma)\big)
$$
and since we can choose the Clifford multiplication by $N$ such that
\begin{eqnarray*}
\gamma^\Omega(N)=-i
\begin{pmatrix}
\Id & 0\\
0 & -\Id
\end{pmatrix}
\end{eqnarray*}
we finally get
\begin{eqnarray*}
\D^\pm=
\begin{pmatrix}
\pm\frac{n}{2}\Id & D^\Sigma  & 0 & 0 \\
D^\Sigma & \mp\frac{n}{2}\Id & 0 & 0 \\
0 & 0 & \mp\frac{n}{2}\Id & D^\Sigma \\
0 & 0 & D^\Sigma & \pm\frac{n}{2}\Id
\end{pmatrix}.
\end{eqnarray*}
This expression clearly shows that these operators only depend  on the Riemannian metric and the spin structure on $\Sigma$. Here the matrix blocks are defined with respect to the decomposition 
$$\Eb\cong\big(\mathbb{S}^+(\Sigma)\oplus\mathbb{S}^-(\Sigma)\big)\oplus\big(\mathbb{S}^+(\Sigma)\oplus\mathbb{S}^-(\Sigma)\big).$$

We summarize the preceding discussion by
\begin{proposition}
The bundle $(\E\Omega,\gamma,\nabla)$ is a Dirac bundle equipped with a chirality operator $\omega$ whose associated Dirac-type operator $D:=\gamma\circ\nabla$ is a first order elliptic differential operator. The restricted triplet $(\Eb,\mult,\nb)$ is also a Dirac bundle for which the spinorial Gauss formula
\begin{eqnarray}\label{S-G-F-E}
\nb_X\Psi=\nabla_X\Psi - \frac{1}{2}\mult(AX)\Psi
\end{eqnarray} 
holds for all $\Psi\in\Gamma(\Eb)$ and $X\in\Gamma(T\Sigma)$ and such that
\begin{eqnarray}\label{D-D-E}
\D\Psi=\frac{n}{2}H\Psi-\gamma(N)D\Psi-\nabla_N\Psi
\end{eqnarray}
and 
\begin{eqnarray}\label{DE-C}
\D\big(\gamma(N)\Psi\big)=-\gamma(N)\D\Psi
\end{eqnarray}
where $\D:=\mult\circ\nb$ is the extrinsic Dirac-type operator on $\Eb$. Moreover, the Dirac-type operators $\D^\pm:=\D\pm\frac{n}{2}i\gamma(N)\Id_{\Eb}$ are first order differential operators which only depend on the Riemannian and spin structures of $\Sigma$.
\end{proposition}


\subsection{The Hyperbolic Reilly formula}


We first recall the hyperbolic version of the Schr\"odinger-Lichnerowicz formula on the spinor bundle where a proof can be found in \cite{AD}, \cite{HMR2} or \cite{minoo}
\begin{eqnarray*}
\int_{\Omega}\big(\frac{1}{4}\big(R+n(n+1)\big)|\psi|^2-\frac{n}{n+1}|D^{\Omega}_{\pm}\psi|^2\big)\,d\Omega 
\leq\int_{\Sigma} \big(\<\D^\Sigma_\pm\psi,\psi\> -\frac{n}{2}{H}|\psi|^2\big)\,d\Sigma
\end{eqnarray*}
for all $\psi\in\Gamma(\So)$ and where $D^{\Omega}_\pm:=D^\Omega\mp\frac{n+1}{2}i\Id$ and $\D^\Sigma_\pm:=\D^\Sigma\pm\frac{n}{2}i\gamma^\Omega(N)$. Moreover equality occurs if and only if $\psi$ is a twistor-spinor and the scalar curvature of $\Omega$ is constant equal to $-n(n+1)$. Recall that a twistor-spinor on $\So$ is a smooth spinor field such that $P^\Omega_X\psi=0$ for all $X\in\Gamma(T\Omega)$ where the operator $P^\Omega$ is the twistor operator (also called Penrose operator) defined for all $\psi\in\Gamma(\mathbb{S}\Omega)$ by
\begin{eqnarray*}
P^{\Omega}_X\psi:=\nabla_X^\Omega\psi+\frac{1}{n+1}\gamma^\Omega(X)D^\Omega\psi,
\end{eqnarray*}
(for more details, we refer to \cite{BFGK}). We now extend the above Hyperbolic Reilly Inequality to  sections of the Dirac bundle $\E\Omega$. For this, we define the twistor operator on $\E\Omega$ by
$$
P_X:=\nabla_X+\frac{1}{n+1}\gamma(X)D=\left\lbrace
\begin{array}{ll}
P_X^\Omega \qquad& {\rm if }\,\, n+1=2m\\ \\
P_X^\Omega \oplus P_X^\Omega \qquad& {\rm if }\,\, n+1=2m+1
\end{array}
\right.
$$
and a section $\Psi\in\Gamma(\E\Omega)$ such that $P_X\Psi=0$ for all $X\in\Gamma(T\Omega)$ will be called a twistor-spinor on $\E\Omega$. Then it is a simple exercise to check that the following formula holds on $\E\Omega$:
\begin{proposition}
Let $\Omega$ be a compact and connected $(n+1)$-dimensional Riemannian spin manifold with boundary $\Sigma$. Assume that the scalar curvature of $\Omega$ satisfies $R\geq-n(n+1)$, then for all $\Psi\in\Gamma(\E\Omega)$, we have
\begin{eqnarray}\label{h-s-r}
-\frac{n}{n+1}\int_{\Omega}|D^\pm\Psi|^2\,d\Omega 
\leq\int_{\Sigma} \big(\<\D^\pm\Psi,\Psi\> -\frac{n}{2}H|\Psi|^2\big)\,d\Sigma.
\end{eqnarray}
Moreover equality occurs if and only if $\Psi$ is a twistor-spinor on $\E\Omega$ and $R=-n(n+1)$. Here
$\D^\pm$ are defined in $(\ref{MEDTO})$ and $D^\pm$ are the modified Dirac-type operators acting on sections of $\E\Omega$ defined by:
\begin{eqnarray}\label{MDTO}
D^\pm:=D\mp\frac{n+1}{2}i\,\Id_{\E\Omega}.
\end{eqnarray} 
\end{proposition}


\subsection{A boundary-value value problem for the Dirac-type operator $D^+$}\label{an-set}


In this section, we introduce the boundary condition which we will need and prove its ellipticity for a Dirac-type operator acting on  $\Gamma(\E\Omega)$. It turns out that this condition is well-known for even dimensional Riemannian spin manifolds: this is the condition associated with a chirality operator (see \cite{HMR1} for example). 
Here we extend it for odd dimensional Riemannian spin manifolds. Note that, as explained in the previous section, we are not working on the  spinor bundle $\So$ since this boundary condition does not yield to an elliptic boundary condition for the fundamental Dirac operator $D^\Omega$ on $\Omega$. 

\vspace{0.2cm}

Since the modified Dirac-type operators $D^\pm$  (see (\ref{MDTO})) acting on sections of $\E\Omega$   are zero order deformations of the Dirac operator, they define first order elliptic differential operators whose $L^2$-formal adjoints are  $(D^\pm)^*=D^\mp$. This last fact is an obvious consequence of the following integration by parts formula
\begin{eqnarray}\label{s-f}
\int_{\Omega}\<D\Psi, \Phi\>\,d\Omega=\int_{\Omega}\<\Psi,D\Phi\>\,d\Omega-\int_\Sigma\<\gamma(N)\Psi,\Phi\>\,d\Sigma
\end{eqnarray}
for all $\Psi,\Phi\in\Gamma(\E\Omega)$ and where $d\Omega$ (resp. $d\Sigma$) is the Riemannian volume element of $\Omega$ (resp. $\Sigma$). It is then easy to see that we are in the standard setup examined by B\"ar and Ballmann (see page 
$5$ of \cite{BB} for a precise definition of this setting). On the other hand, the fiber preserving endomorphism 
$$G=\gamma(N)\omega:\Gamma(\Eb)\rightarrow\Gamma(\Eb),$$
acting on sections of the restricted bundle, is self-adjoint with respect to the pointwise Hermitian scalar product, whose square is the identity. Here $\omega$ is the chirality operator defined in Section \ref{G-S}. The map $G$ has two eigenvalue $\pm 1$ whose corresponding eigenspaces are interchanged by the isomorphism $\gamma(N)$.  Then we consider the two non trivial eigensubbundles $\mathcal V^\pm$ over $\Sigma$ corresponding to the $\pm 1-$eigenvalues  of the map $G$ so that the following decomposition holds$$
\Eb=\mathcal V^+\oplus \mathcal V^-.
$$
The pointwise projections on $\mathcal V^\pm$ are given by
\begin{equation}\label{p-c-o}
\begin{array}{r@{}ccl}
 P_\pm \colon  & L^2(\E\Omega) &\lto & L^2(\mathcal V^\pm) \\[.05cm]  
  & \Psi &\longmapsto  & P_\pm\Psi :=\frac{1}{2}(\Id\pm\gamma(N)\omega)\Psi, 
\end{array}
\end{equation}
where $L^2(\E\Omega)$ (resp. $L^2(\mathcal V^\pm)$) denotes the space of $L^2$-integrable sections of $\E\Omega$ (resp. $\mathcal V^\pm$). Using the properties (\ref{c-o-1}) and (\ref{c-o-2}) of $\omega$, we easily see that for $X\in\Gamma(T\Sigma)$ and $\Psi\in\Gamma(\Eb)$, we have $G(\mult(X)\Psi)=\mp\mult(X)G\Psi$ and then $\mult(X)$ interchanges $\mathcal V^+$ and $\mathcal V^-$. So from Corollary $7.23$ and Proposition $7.24$ in \cite{BB}, we have
\begin{proposition}
The pointwise orthogonal decomposition $\Eb=\mathcal V^+\oplus \mathcal V^-$ induces local boundary conditions for $D^+$. In particular, the operator
\begin{eqnarray*}
D^+:{\rm Dom}(D^+):=\{\Psi\in H^2_1\,: \,P_\pm(\Psi_{|\Sigma})=0\}\longrightarrow L^2(\E\Omega)
\end{eqnarray*}
is Fredholm and if $\Phi$ is a smooth section of $\E\Omega$, then any $H^2_1$-solutions of 
$$
\left\lbrace
\begin{array}{ll}
D^+\Psi=\Phi \quad& \,\text{on}\;\Omega\\
P_\pm\Psi_{|\Sigma}=0 \quad& \,\text{along}\;\Sigma,
\end{array}
\right.
$$
is smooth up to the boundary. Here $H^2_1$ stands for the Sobolev space of $L^2$-spinors with weak $L^2$-covariant derivatives.
\end{proposition}
It is clear that the same result holds for the Dirac-type operator $D^-$. Next we only consider the operator $D^+$ since it is straightforward to check that all the following results also hold for $D^-$.

Now we want to prove that the Dirac operator $D^+$ defines an isomorphism between the space 
$${\rm Dom}_\pm(D^+):=\{\Psi\in H^2_1(\E\Omega)\, : \,P_\pm\Psi_{|\Sigma}=0\}$$ 
onto $L^2(\E\Omega)$, where $P_\pm$ is the projection given by (\ref{p-c-o}). We now denote by $D^+_\pm$ the Dirac-type operator defined on the domain ${\rm Dom}_\pm(D^+)$. We have
\begin{proposition}\label{dir-iso}
Let $\Omega$ be a compact domain with smooth boundary in a $(n+1)$-dimensional Riemannian spin manifold. The Dirac-type operator $D^+$ with domain ${\rm Dom}_\pm(D^+)$ is an isomorphism onto the space of square integrable sections of $\E\Omega$. In particular, for all $\Phi\in\Gamma(\E\Omega)$, there exists a unique smooth section $\Psi\in\Gamma(\E\Omega)$ such that
\begin{equation}\label{iso-chi}
\left\lbrace
\begin{array}{ll}
D^+\Psi=\Phi \quad& \,\text{on}\;\Omega\\
P_\pm\Psi_{|\Sigma}=0 \quad& \,\text{along}\;\Sigma.
\end{array}
\right.
\end{equation}
\end{proposition}

\pf
From the Stokes' formula (\ref{s-f}) we have for all $\Psi,\Phi\in\Gamma(\E\Omega)$
\begin{eqnarray*}
\int_\Omega\<D^+\Psi,\Phi\>\,d\Omega=\int_\Omega\<\Psi,D^-\Phi\>\,d\Omega-\int_\Sigma\<\gamma(N)\Psi,\Phi\>\,d\Sigma.
\end{eqnarray*}
On the other hand for all $\Psi\in\Gamma(\E\Omega)$ 
\begin{eqnarray*}
P_\pm\Psi_{|\Sigma}=0\Longleftrightarrow P_\mp\big(\gamma(N)\Psi_{|\Sigma}\big)=0,
\end{eqnarray*} 
then the boundary term of the previous identity vanishes for all $\Psi \in {\rm Dom}_\pm(D^+)$, hence 
$(D^+_\pm)^*=D^-_\pm$. Since 
\begin{eqnarray}\label{coker}
{\rm CoKer}(D^+_\pm)\simeq{\rm Ker}(D^+_\pm)^*\simeq{\rm Ker}(D^-_\pm),
\end{eqnarray}
we only have to show that ${\rm Ker}(D^+_\pm)$ and ${\rm Ker}(D^-_\pm)$ are reduced to zero to conclude that $D^+_\pm$ is an isomorphism. So if $\Psi\in\Gamma(\E\Omega)$ is in the kernel of $D^+_\pm$ that is
$$
\left\lbrace
\begin{array}{ll}
D^+\Psi=0 \quad& \,\text{on}\;\Omega\\
P_\pm\Psi_{|\Sigma}=0 \quad& \,\text{along}\;\Sigma,
\end{array}
\right.
$$
then, from the ellipticity of the boundary condition $P_\pm$, it has to be smooth up to the boundary. Moreover since $D^+\Psi=0$ we have on one hand
\begin{eqnarray*}
\int_\Omega\<D\Psi,\Psi\>\,d\Omega= i\frac{n+1}{2}\int_\Omega|\Psi|^2\,d\Omega,
\end{eqnarray*} 
and on the other hand, an integration by parts leads to
\begin{eqnarray*}
\int_\Omega\<D\Psi,\Psi\>\,d\Omega=\int_\Omega\<\Psi,D\Psi\>\,d\Omega=-i\frac{n+1}{2}\int_\Omega|\Psi|^2\,d\Omega. 
\end{eqnarray*}
In other words, we showed that
\begin{eqnarray*}
(n+1)i\int_\Omega|\Psi|^2\,d\Omega=0
\end{eqnarray*}
which implies that $\Psi\equiv 0$ on $\Omega$. We conclude that the kernel of $D^+_\pm$ is trivial and by using a similar argument, we also get that ${\rm Ker}(D^-_\pm) =\{0\}.$  Then from (\ref{coker}), the operator $D^+_\pm$ is an isomorphism. From this fact, it is obvious to see that for all $\Phi\in\Gamma(\E\Omega)$, there exists a unique smooth solution of (\ref{iso-chi}).

\qed

As a consequence, we prove that the associated non-homogeneous boundary-value problem has a unique smooth solution: 
\begin{corollary}\label{H-Ext}
Let $\Sigma$ be a hypersurface bounding a compact domain $\Omega$ in an $(n+1)$-dimensional Riemannian spin manifold. Then for all $\Phi\in\Gamma(\Eb)$, 
 there exists a non trivial smooth section $\Psi\in\Gamma(\E\Omega)$, solution of the boundary-value problem
\begin{equation}\label{harm-ext}
\left\lbrace
\begin{array}{ll}
D^+\Psi=0 & \,\text{on}\;\Omega\\
P_+\Psi_{|\Sigma}=P_+\Phi & \,\text{along}\;\Sigma.
\end{array}
\right.
\end{equation}
The same conclusion holds for the boundary condition $P_-$.
\end{corollary}

\pf 
 Let $\widehat{\Phi}$ be a smooth extension of $\Phi$ on $\Omega$. From Proposition \ref{dir-iso}, there exists a smooth solution $\widehat{\Psi}\in\Gamma(\E\Omega)$ to the boundary-value problem
$$
\left\lbrace
\begin{array}{ll}
D^+\widehat{\Psi}=-D^+\widehat{\Phi} \quad& \,\text{on}\;\Omega\\
P_\pm\widehat{\Psi}_{|\Sigma}= 0 \quad& \,\text{along}\;\Sigma.
\end{array}
\right.
$$
It is then straightforward to see that $\Psi:=\widehat{\Psi}+\widehat{\Phi}$ is a smooth section of $\E\Omega$ which satisfies (\ref{harm-ext}).
\qed


\section{The holographic principle for Dirac bundles}


In this section, we prove a holographic principle for the existence of an imaginary Killing spinor. This result is the hyperbolic counterpart of a similar principle for parallel spinor fields proved by the first two authors in \cite{HM1}.
\begin{theorem}\label{holo-im}
Let $\Omega$ be a compact, connected Riemannian spin manifold with smooth boundary $\Sigma$. Assume that the scalar curvature of $\Omega$ satisfies $R\geq-n(n+1)k^2$ for some $k>0$ and the mean curvature $H$ of $\Sigma$ is positive, then
for all $\Phi\in\Gamma(\Eb)$, one has
\begin{eqnarray}\label{holoim}
\int_\Sigma\Big(\frac{1}{H}|\D^+\Phi|^2-\frac{n^2}{4}H|\Phi|^2\Big)\,d\Sigma\geq 0 .
\end{eqnarray}
 Moreover, equality occurs for $\Phi\in\Gamma(\Eb)$ if and only if there exists two imaginary Killing spinor fields $\Psi^+$, $\Psi^-\in\Gamma(\Eb)$ with Killing number $-(i/2)$ such that $P_+\Psi^+=P_+\Phi$ and $P_-\Psi^-=P_-\Phi$.
\end{theorem}

\begin{remark}\label{Class-Im}
{\rm In the previous result, a smooth section $\Phi^\pm\in\Gamma(\E\Omega)$ is called an imaginary Killing spinor on $\E\Omega$ with Killing number $\pm(i/2)$ if it satisfies the equation 
\begin{eqnarray*}
\nabla_X\Phi^\pm=\pm\frac{i}{2}\gamma(X)\Phi^\pm
\end{eqnarray*}
for all $X\in\Gamma(T\Omega)$. It is clear that if $\Omega$ is an even dimensional manifold, the existence of an imaginary Killing spinor on $\E\Omega$ is equivalent to the existence of an imaginary  Killing spinor on $\mathbb{S}\Omega$ since in this case,  $\E \Omega=\mathbb{S}\Omega$. If the dimension of $\Omega$ is odd, the existence of one imaginary Killing spinor with Killing number $\pm(i/2)$ is enough to ensure that $\E\Omega$ carries two imaginary Killing spinors with Killing number $(i/2)$ {\rm and} $-(i/2)$. Indeed, it is immediate to check that if $\phi$ denotes such a spinor field on $\mathbb{S}\Omega$, then the fields defined on $\E\Omega$ by $\Phi^+=(\phi,0)$ and $\Phi^-=(0,\phi)$ are imaginary Killing spinors on $\E\Omega$ whose Killing number are respectively $\pm(i/2)$ and $\mp(i/2)$. Moreover, they satisfy $|\Phi^\pm|^2=|\phi|^2$ and they have no zero since imaginary Killing spinors in $\So$ have no zero (see \cite{B1} or \cite{B2} for example). }
\end{remark}

The choice of the boundary condition heavily relies on its behavior with respect to the modified Dirac-type operator $\D^\pm$. We first state the main properties needed here to prove our main result. 
\begin{lemma}
The Dirac-type operator $\D^\pm$ defined for all $\Phi\in\Gamma(\Eb)$ by:\begin{eqnarray*}
\D^\pm\Phi:=\D\Phi\pm\frac{n}{2}i\gamma(N)\Phi
\end{eqnarray*}
are first order elliptic differential operators which are self-adjoint with respect to the $L^2$-scalar product on $\Eb$. Moreover for all $\Phi\in\Gamma(\Eb)$, we have:
\begin{eqnarray}\label{D-P-C}
\D^+(P_\pm\Phi)=P_\mp(\D^+\Phi)
\end{eqnarray}
and so in particular:
\begin{eqnarray}\label{bound-exp}
\int_\Sigma\<\D^+\Phi,\Phi\>\,d\Sigma=2\int_{\Sigma}{\rm Re}\<\D^+(P_+\Phi),P_-\Phi\>\,d\Sigma.
\end{eqnarray}
\end{lemma}

\pf
First note that, since $\D^+$ is a zero order deformation of the first order elliptic differential operator $\D$, it is also a first order elliptic operator. 
Then note that the endomorphism $i\gamma(N)$ of $\Eb$ is symmetric with respect to the pointwise Hermitian scalar product $\<\,,\,\>$, so that we easily compute for all $\Phi_1$, $\Phi_2\in\Gamma(\Eb)$:
\begin{eqnarray*}
\int_\Sigma\<\D^+\Phi_1,\Phi_2\>\,d\Sigma=\int_\Sigma\<\D\Phi_1+\frac{n}{2}i\gamma(N)\Phi_1,\Phi_2\>\,d\Sigma=\int_\Sigma\<\Phi_1,\D^+\Phi_2\>\,d\Sigma
\end{eqnarray*} 
since $\D$ is $L^2$-self-adjoint. This proves the first assertion. A straightforward computation shows that $\gamma(N)P_\pm=P_\mp\gamma(N)$ and then the skew-commutativity rule (\ref{DE-C}) gives (\ref{D-P-C}). Now every section of $\Eb$ can be decomposed into $\Phi=P_+\Phi+P_-\Phi$ and since this decomposition is pointwise orthogonal, we compute using (\ref{D-P-C}):
\begin{eqnarray*}
\int_\Sigma\<\D^+\Phi,\Phi\>\,d\Sigma & =
& \int_\Sigma \<P_-(\D^+\Phi),P_-\Phi\>\,d\Sigma+\int_\Sigma\<P_+(\D^+\Phi),P_+\Phi\>\,d\Sigma\\
& = & \int_\Sigma\<\D^+(P_+\Phi),P_-\Phi\>\,d\Sigma+\int_\Sigma\<P_-\Phi,\D^+(P_+\Phi)\>\,d\Sigma\\
& = & 2\int_{\Sigma}{\rm Re}\<\D^+(P_+\Phi),P_-\Phi\>\,d\Sigma.
\end{eqnarray*}

\qed

\begin{proposition}\label{holo1}
Let $\Omega$ be a compact spin Riemannian manifold with scalar curvature $R\geq -n(n+1)$ and whose boundary $\Sigma$ has positive mean curvature $H$. For any section $\Phi$ of the restricted Dirac bundle $\Eb$, one has  
\begin{eqnarray}\label{holo1+}
0\leq\int_{\Sigma}\big(\frac{1}{H}|\D^+ P_+\Phi|^2-\frac{n^2}{4}H|P_+\Phi|^2\big)\,d\Sigma.
\end{eqnarray}
Moreover, equality holds if and only if there exists an imaginary Killing spinor $\Psi^+\in\Gamma (\E\Omega )$ such that $P_+\Psi^+=P_+\Phi$ along the boundary.
\end{proposition}

\pf
Take any spinor field $\Phi\in\Gamma(\Eb)$ on the hypersurface and consider the following boundary-value problem
\begin{equation}\label{bound-val}
\left\{
\begin{array}{lll}
{D}^+ \Psi^+&=0 \qquad&\hbox{ {\rm on} } \Omega \\
P_+\Psi^+_{|\Sigma} &=P_+ \Phi \qquad&\hbox{ {\rm on} }\Sigma 
\end{array}
\right.
\end{equation}
for the Dirac-type operator ${D}^+$ and the boundary condition $P_+$. The existence and uniqueness of a smooth solution $\Psi^+\in \Gamma(\E\Omega)$ for this boundary-value problem is ensured by Corollary \ref{H-Ext}. On the other hand, since we assume that $R\geq -n(n+1)$, we can apply the hyperbolic Reilly inequality (\ref{h-s-r}) to $\Psi^+$ to get the following inequality 
\begin{eqnarray*}
0\leq\int_\Sigma\big(\<\D^+\Psi^+,\Psi^+\>-\frac{n}{2}H|\Psi^+|^2\big)\,d\Sigma.
\end{eqnarray*}
This inequality combined with (\ref{bound-exp}), imply  
\begin{equation}\label{final-ineq-P}
0\le \int_\Sigma \Big(2{\rm Re}\<\D^+P_+\Psi^+,P_-\Psi^+\>
-\frac{n}{2}H|P_+\Psi^+|^2-\frac{n}{2}H|P_-\Psi^+|^2\Big)\,d\Sigma.
\end{equation}
Since we assume that the mean curvature $H>0$, we can write
\begin{eqnarray*}
&{\displaystyle
0\le \big|\sqrt{\frac{2}{nH}}\;\D^+ P_+\Psi^+-\sqrt{\frac{nH}{2}} \;P_-\Psi^+\big|^2 
= } \nonumber \\
&{\displaystyle \frac{2}{nH}|\D^+ P_+\Psi^+|^2+\frac{nH}{2}|P_-\Psi^+|^2-2 {\rm Re} \langle\D P_+\Psi^+,P_-\Psi^+\rangle.}\nonumber
\end{eqnarray*}
In other words, we have
\beQ
2 {\rm Re} \langle\D^+ P_+\Psi^+, P_-\Psi^+\rangle -\frac{nH}{2}|P_-\Psi^+|^2 \le \frac{2}{nH}|\D^+ P_+\Psi^+|^2,
\eeQ
which, when combined with Inequality (\ref{final-ineq-P}), knowing that $P_+\Psi^+=P_+\Phi$, imply Inequality (\ref{holo1+}).  

Assume now that equality is achieved, then the spinor field $\Psi^+\in\Gamma(\E\Omega)$ which satisfies the boundary-value problem (\ref{bound-val}) is in fact a twistor-spinor since we have equality in the hyperbolic Reilly formula (\ref{h-s-r}). Moreover, since the condition $D^+\Psi^+=0$ translates to $D\Psi^+=\frac{n+1}{2}i\Psi^+$, the section $\Psi^+$ is in fact an imaginary Killing spinor on $\E\Omega$ with Killing number $-(i/2)$. Moreover, it is obvious that $P_+\Psi^+_{|\Sigma}=P_+\Phi$ as asserted. 

Conversely, if $\Psi^+$ is an imaginary Killing spinor on $\E\Omega$ then from 
 (\ref{D-D-E}) we compute
\begin{eqnarray*}
\D\Psi^+ & = & \frac{n}{2}H\Psi^+-\gamma(N)D\Psi^+-\nabla_N\Psi^+\\
& = & \frac{n}{2}H\Psi^+-\frac{n}{2}i\gamma(N)\Psi^+
\end{eqnarray*}
which can be written as $\D^+\Psi^+=\frac{n}{2}H\Psi^+$. Now we decompose the section $\Psi^+$ with respect to $P_+$ and $P_-$ and thus the relation (\ref{D-P-C}) yields
\begin{eqnarray}\label{s-d-e}
\D^+(P_\pm\Psi^+)=\frac{n}{2}HP_\mp\Psi^+.
\end{eqnarray}
Moreover, from the $L^2$-self-adjointness of $\D^+$ and (\ref{s-d-e}), we get
\begin{eqnarray*}
\frac{n}{2}\int_\Sigma H|P_-\Psi^+|^2\, d\Sigma & = & \int_\Sigma\<\D^+P_+\Psi^+,P_-\Psi^+\>\,d\Sigma\\
& = & \int_\Sigma\<P_+\Psi^+,\D^+P_-\Psi^+\>\,d\Sigma\\
& = & \frac{n}{2}\int_\Sigma H|P_+\Psi^+|^2\, d\Sigma
\end{eqnarray*}
that is
\begin{eqnarray}\label{est-bounb}
\int_\Sigma H|P_-\Psi^+|^2\, d\Sigma = \int_\Sigma H|P_+\Psi^+|^2\, d\Sigma.
\end{eqnarray}
Finally, using (\ref{s-d-e}) and (\ref{est-bounb}), it follows
\begin{eqnarray*}
\int_{\Sigma}\big(\frac{1}{H}|\D P_+\Psi^+|^2-\frac{n^2}{4}H|P_+\Psi^+|^2\big)\,d\Sigma=\frac{n^2}{4}\int_\Sigma H(|P_-\Psi^+|^2-|P_+\Psi^+|^2)\,d\Sigma
\end{eqnarray*}
so that equality is achieved in (\ref{holo1+}).

\qed 

We can mimic this proof step by step to get the counterpart of this result for the boundary condition $P_-$
\begin{proposition}\label{holo2}
Let $\Omega$ be a compact spin Riemannian manifold with scalar curvature $R\geq -n(n+1)$, whose boundary $\Sigma$ has positive mean curvature $H$. For any section $\Phi$ of the restricted Dirac bundle $\Eb$, one has  
\begin{eqnarray}\label{holo2+}
0\leq\int_{\Sigma}\big(\frac{1}{H}|\D^+ P_-\Phi|^2-\frac{n^2}{4}H|P_-\Phi|^2\big)\,d\Sigma.
\end{eqnarray}
Moreover, equality holds if and only if, there exists an imaginary Killing spinor $\Psi^-$ on $\E\Omega$ such that $P_-\Psi^-=P_-\Phi$, along the boundary.
\end{proposition}

\noindent {\textit {Proof of {\bf Theorem \ref{holo-im}}:} }  
By Propositions \ref{holo1} and \ref{holo2}, 
the  field $\Phi\in\Gamma(\Eb)$ satisfies  inequalities (\ref{holo1+}) and (\ref{holo2+}). Summing these estimates and using the relation (\ref{D-P-C}) gives the result. The equality case also follows directly from the characterization of the equality cases in Propositions \ref{holo1} and \ref{holo2}. 

\qed

Now making use of the restriction to the hypersurface of an imaginary Killing spinor field, we get
\begin{theorem}\label{OneIm}
Let $(\Omega^{n+1},g)$ be a compact, connected $(n+1)$-dimensional Riemannian spin manifold with smooth boundary $\Sigma$. Assume that the scalar curvature of $\Omega$ satisfies $R\geq -n(n+1)$ and that the mean curvature $H$ of $\Sigma$ is positive.  Suppose furthermore that $\Sigma$ admits an isometric and isospin immersion $F$ into another $(n+1)$-dimensional Riemannian spin manifold $(\Omega_0,g_0)$ endowed with a non trivial $\pm(i/2)$-imaginary Killing spinor field $\Phi^\pm\in\Gamma(\E \Omega_0)$ and denote by $H_0$ the mean curvature of this immersion. Then the following inequality holds
\begin{eqnarray}\label{E-M-V1}
\int_\Sigma\Big(\frac{H_0^2-H^2}{H}\Big)|\Phi^\pm|^2\,d\Sigma\geq 0
\end{eqnarray}
and equality occurs if and only if both immersions have the same shape operators and $\Sigma$ is connected.
\end{theorem}

\noindent {\textit {Proof of {\bf Theorem \ref{OneIm}} :} } We only consider the case where $\Phi^-\in\Gamma(\E \Omega_0)$ is an imaginary Killing spinor with Killing number $-(i/2)$. If $\Sigma_0$ is a connected component of the boundary $\Sigma$, then, by taking the restriction of the imaginary Killing spinor $\Phi^-\in\Gamma(\E \Omega_0)$ to $\Sigma_0$, we get the existence of a section $\Phi_0^-:=\Phi^-_{|\Sigma_0}$ which satisfies the intrinsic Dirac-type equation
\begin{eqnarray}\label{DirEq}
\D^+\Phi^-_0=\frac{n}{2}H_0\Phi^-_0.
\end{eqnarray} 
Now we extend the section $\Phi^-_0$ on $\Sigma$ in such a way that its extension, also denoted by $\Phi_0^-\in\Gamma(\Eb)$, vanishes on $\Sigma-\Sigma_0$. Then putting this spinor field into  (\ref{holoim}) gives the estimate (\ref{E-M-V1}). Assume now that equality is achieved, then from the equality case of (\ref{holoim}), there exists two imaginary Killing spinor fields $\Psi^+$, $\Psi^-\in\Gamma(\E\Omega)$ with Killing number $-(i/2)$ such that $P_+\Psi^+=P_+\Phi_0^-$ and $P_-\Psi^-=P_-\Phi_0^-$. Using (\ref{DirEq}), (\ref{D-P-C}) and Formula (\ref{D-D-E}), we have
\begin{eqnarray}\label{id1}
H_0P_+\Phi_0^-=\frac{2}{n}\D^+(P_-\Phi_0^-)=\frac{2}{n}\D^+(P_-\Psi^-)=HP_+\Psi^-.
\end{eqnarray}
Similarly, we obtain
\begin{eqnarray}\label{id2}
H_0P_-\Phi_0^-=\frac{2}{n}\D^+(P_+\Phi_0^-)=\frac{2}{n}\D^+(P_+\Psi^+)=HP_-\Psi^+.
\end{eqnarray}
Applying the operator $\D$ to the first and last terms of (\ref{id1}), we get
\begin{eqnarray*}
\mult(\nabla^\Sigma H_0)P_+\Phi_0^-+\frac{n}{2}H^2_0P_-\Phi_0^-=\mult(\nabla^\Sigma H)P_+\Psi^-+\frac{n}{2}H^2P_-\Psi^-
\end{eqnarray*}
which, using again the equalities above, finally gives
\begin{eqnarray*}
\mult(\nabla^\Sigma H_0)P_+\Phi_0^-+\frac{n}{2}H^2_0P_-\Phi_0^-=\frac{H_0}{H}\mult(\nabla^\Sigma H)P_+\Phi_0^-+\frac{n}{2}H^2P_-\Phi_0^-.
\end{eqnarray*}
The same argument applied to (\ref{id2}) yields
\begin{eqnarray*}
\mult(\nabla^\Sigma H_0)P_-\Phi_0^-+\frac{n}{2}H^2_0P_+\Phi_0^-=\frac{H_0}{H}\mult(\nabla^\Sigma H)P_-\Phi_0^-+\frac{n}{2}H^2P_+\Phi_0^-,
\end{eqnarray*}
so that the sum of the last two formulae implies
\begin{eqnarray*}
\mult(\nabla^\Sigma H_0)\Phi_0^-+\frac{n}{2}H^2_0\Phi_0^-=\frac{H_0}{H}\mult(\nabla^\Sigma H)\Phi_0^-+\frac{n}{2}H^2\Phi_0^-.
\end{eqnarray*}
Moreover, since the spinor fields $\mult(\nabla^\Sigma H_0)\Phi_0^-$ and $\mult(\nabla^\Sigma H)\Phi_0^-$ are both orthogonal to $\Phi_0^-$, and since the spinor $\Phi_0^-$ has no zeros on $\Sigma_0$ (see Remark \ref{Class-Im}), we deduce that $H_0^2=H^2$ and $H\nabla^\Sigma H_0=H_0\nabla^\Sigma H$. From these facts, we conclude that $H_0$ has no zeros since $H$ is positive and so we may assume that $H_0=H$. Using this equality in (\ref{id1}) and (\ref{id2}) gives $\Phi_{0|\Sigma}^-=\Psi^+_{|\Sigma}=\Psi^-_{|\Sigma}$. By definition, we have $\Phi_{0|\Sigma-\Sigma_0}=0$ on $\Sigma-\Sigma_0$, thus
\begin{eqnarray*}
P_+\Psi^+=P_+\Phi_0^-=0\quad{\rm and}\quad P_-\Psi^-=P_-\Phi^-_0=0.
\end{eqnarray*} 
Applying the operator $\D^+$ to these equalities and using (\ref{D-P-C}) and (\ref{D-D-E}), we get
\begin{eqnarray*}
0=\D^+(P_+\Psi^+)=\frac{n}{2}HP_-\Psi^+,\quad 0=\D^+(P_-\Psi^-)=\frac{n}{2}HP_+\Psi^-
\end{eqnarray*}
and since $H>0$, we deduce 
$$\Psi^+_{|\Sigma-\Sigma_0}=\Psi^-_{|\Sigma-\Sigma_0}=0.$$
However, since $\Psi^+$ and $\Psi^-$ are imaginary Killing spinors on $\E\Omega$, they have no zeros, so this is impossible unless $\Sigma=\Sigma_0$ is connected.\\

Finally, as another consequence of the preceding argument, we have that $\Phi_0^-$ is the restriction to $\Sigma$ of $\Psi^+$ (and $\Psi^-$) via the embedding of $\Sigma$ as the boundary of $\Omega$ and of $\Phi^-$ via the immersion of $\Sigma$ in $\Omega_0$. Then we can apply the spinorial Gauss formula (\ref{S-G-F-E}) for the first immersion, that is 
\begin{eqnarray}\label{Im1}
\nb_X\Psi^+=-\frac{i}{2}\gamma(X)\Psi^+-\frac{1}{2}\mult(AX)\Psi^+
\end{eqnarray}
for all $X\in\Gamma(T\Sigma)$ (here $A$ is the second fundamental form of $\Sigma\hookrightarrow\Omega$), as well as 
\begin{eqnarray}\label{Im2}
\nb_X\Phi_0^-=-\frac{i}{2}\gamma^0(X)\Phi_0^--\frac{1}{2}\mult(A_0X)\Phi_0^-
\end{eqnarray}
for the second immersion $\Sigma\hookrightarrow\ \Omega_0$. The notation $\gamma^0$ stands for the Clifford multiplication on $\E\Omega_0$. Now we claim that
\begin{eqnarray}\label{Im3}
\gamma(X)\Psi^+=\gamma^0(X)\Phi_0^-
\end{eqnarray}
for all $X\in\Gamma(T\Sigma)$. Indeed from Section \ref{Sp-Ch}, we have seen that we can choose $\gamma^0(N_0)$ and $\gamma(N)$ such that $\gamma^0(N_0)\Phi=\gamma(N)\Phi$ for all $\Phi\in\Gamma(\Eb)$ and thus for all $X\in\Gamma(T\Sigma)$,
\begin{eqnarray*}
\gamma(X)\Phi=-\mult(X)\gamma(N)\Phi=-\mult(X)\gamma^0(N_0)\Phi=\gamma^0(X)\Phi.
\end{eqnarray*}
Using the fact that $\Psi^+_{|\Sigma}=\Phi_{0|\Sigma}^-$ in (\ref{Im1}) and (\ref{Im2}) with the relation (\ref{Im3}) finally give
\begin{eqnarray*}
\mult(A_0X-AX)\Phi_{0|\Sigma}^-=0
\end{eqnarray*}
for $X$ tangent to $\Sigma$, and since $\Phi_0^-$ has no zeros, we get $A_0=A$.\\

The converse is clear. If the two shape operators $A$ and $A_0$ coincide, then the corresponding traces $nH$ and $nH_0$ taken with respect to the common induced metric should be equal. Then we have equality in (\ref{E-M-V1}).

\qed


\section{A new quasi-local mass}


We propose here a local version of the positive mass theorem obtained by Wang \cite{Wa1} and Chru{\'s}ciel-Herzlich \cite{CH} for asymptotically hyperbolic manifolds. 


\subsection{The Hyperbolic space and Hypersurfaces}


In this section, we recall some well known facts regarding imaginary Killing spinors of the hyperbolic space ${\mathbb{H}}^{n+1}$. The classification of complete manifolds carrying an imaginary Killing spinor has been obtained by H. Baum in \cite{B1,B2}  (see Remark \ref{Imag-Kil} below). A standard model of the hyperbolic space is the unit ball $\UB$ endowed with the Riemannian metric $g_{{\mathbb{H}}}=f^2g_{\mathbb{E}}$ where $g_{{\mathbb{E}}}$ is the Euclidean metric and $f(x)=2/(1-|x|_{\mathbb{E}}^2)$ for $x\in\UB$. Here $|\,.\,|_{\mathbb{E}}$ denotes the Euclidean norm associated to $g_{\mathbb{E}}$. Since the Riemannian metrics $g_{{\mathbb{H}}}$ and $g_{{\mathbb{E}}}$ are conformally related, we can canonically identify the corresponding spinor bundles $\mathbb{S}{\mathbb{H}}$ and $\mathbb{S}{\mathbb{B}}$. 
Now we consider the $\mathbb{C}^N$-valued constant function on $\UB$  equal to $a\in\mathbb{C}^N$, with $N=2^{[\frac{n+1}{2}]}$, which allows to define a spinor field on the unit Euclidean ball, by setting
\begin{eqnarray}\label{im-kil}
\psi_a^{\pm}(x) :=f^{\frac{1}{2}}(x)(\Id\pm i\gamma^{{\mathbb{E}}}(x))a.
\end{eqnarray}
The spinor field $\psi_a^{\pm}$ induces on ${\mathbb{H}}^{n+1}$ an imaginary Killing spinor field also denoted by $\psi_a^\pm$. In fact, every imaginary Killing spinor on the hyperbolic space can be obtained in such a way. 
\begin{remark}\label{Imag-Kil}
{\rm It is a well-known fact that, after suitable rescaling of the metric, an $(n+1)$-dimensional manifold $P$ with an imaginary Killing spinor has to be Einstein with Ricci curvature $-n$. If $P$ is complete, H. Baum proved in \cite{B1,B2} that it has to be a warped product $\mathbb{R}\times_{\exp}P_0$, i.e. the manifold $\mathbb{R}\times P_0$ is endowed with the metric 
\begin{eqnarray*}
g:=dt^2\oplus e^{2t}g_{P_0}
\end{eqnarray*}  
where $(P_0,g_{P_0})$ is an $n$-dimensional complete Riemannian spin manifold admitting a non-trivial parallel spinor. In case $P_0$ is the Euclidean space $\mathbb{R}^n$, then $P$ is nothing but the hyperbolic space with constant curvature $-1$.}
\end{remark}

In the following, $\Sigma$ is a smooth oriented hypersurface in ${\mathbb{H}}^{n+1}$ whose Weingarten map is denoted by $A_0$, i.e., $A_0(X)=-\nabla_X^\mathbb{H}N_0$ for all $X\in\Gamma(T\Sigma)$, here $\nabla^\mathbb{H}$ is the Levi-Civita connection on $\mathbb{H}^{n+1}$ and $N_0$ is the associated unit inward normal. Now we discuss the existence of {\it imaginary Killing spinors} on $\E\mathbb{H}$ as defined in Remark \ref{Class-Im} and its consequences. For $n+1$ even, the bundle $\E\mathbb{H}$ corresponds to the spinor bundle over the hyperbolic space and this situation is well-known. For the sake of completeness, we include a brief discussion of this case.

\vspace{0.2cm}

\noindent{\it The even dimensional case}

In this case, the bundle $\E\mathbb{H}$ is simply the standard spinor bundle on ${\mathbb{H}}^{2m}$ with Clifford multiplication $\gamma^0=\gamma^{\mathbb{H}}$ and Levi-Civita connection $\nabla^0=\nabla^\mathbb{H}$. Then by the spin Gauss formula (\ref{S-G-F-E}),  the restriction of an imaginary Killing spinor $\psi^\pm_a\in\Gamma(\E{\mathbb{H}})$ to $\Sigma$ satisfies:
\begin{eqnarray*}
\nb_X\psi_a^\pm=\pm\frac{i}{2}\gamma^0(X)\psi_a^\pm-\frac{1}{2}\mult(A_0X)\psi^\pm_a.
\end{eqnarray*}  
In particular, $\psi^\pm_a$ is a solution of the Dirac-type equation
\begin{eqnarray*}
\D\psi^\pm_a=\pm\frac{i}{2}\gamma^0(N_0)\psi_a^\pm+\frac{n}{2}H_0\psi_a^\pm
\end{eqnarray*}
which, by the discussion in Section \ref{Sp-Ch}, translates in an intrinsic way to $\Sigma$ by
\begin{eqnarray*}
\D^\mp\psi^\pm_a=\frac{n}{2}H_0\psi_a^\pm.
\end{eqnarray*}

\noindent{\it The odd dimensional case}

If we assume now that $n=2m$, the vector bundle $\E\mathbb{H}$ is simply two copies of the spinor bundle $\mathbb{S}\mathbb{H}$ with Clifford multiplication $\gamma^0=\gamma^{\mathbb{H}}\oplus-\gamma^{\mathbb{H}}$ and spin Levi-Civita connections $\nabla^0=\nabla^{\mathbb{H}}\oplus\nabla^\mathbb{H}$. In this situation, from the discussion in Remark \ref{Class-Im}, the spinor field defined by $\Psi^\pm_{a}:=(\psi^\pm_a,0)\in\Gamma(\E\mathbb{H})$, where $\psi_a^\pm\in\Gamma(\mathbb{S}\mathbb{H})$ is given by (\ref{im-kil}), satisfies
$$\nabla_X\Psi_{a}^\pm=\pm\frac{i}{2}\gamma(X)\Psi_{a}^\pm$$
for all $X\in\Gamma(T\mathbb{H})$ so that it is an imaginary Killing spinor field on $\E\mathbb{H}$ which, in addition, satisfies $|\Psi^\pm_a|^2=|\psi^\pm_a|^2$. In fact, the spinor field $\Psi^\pm_a$ is characterized by 
\begin{eqnarray*}
\Psi_a^\pm=f^{\frac{1}{2}}(x)(\Id\pm i\gamma_e(x))a
\end{eqnarray*}
where $\gamma_e:=\gamma^{\mathbb{E}}\oplus-\gamma^{\mathbb{E}}$ is the Clifford multiplication on the trivial Dirac bundle $\E\mathbb{B}$ and $a\in\mathbb{C}^N$ is identified with $(a,0)\in\mathbb{C}^N\oplus\mathbb{C}^N$. Now recall from Section \ref{Sp-Ch} that the restricted bundle $\Eb:=\E\mathbb{H}_{|\Sigma}$ can be identified with $\sm\Sigma\oplus\sm\Sigma$, Clifford multiplication $\mult=\gamma^\Sigma\oplus\gamma^\Sigma$ and spin Levi-Civita connection $\nb=\nabla^\Sigma\oplus\nabla^\Sigma$. From these identifications, it is straightforward, using the spinorial Gauss formula (\ref{S-G-F-E}) and the definition of the Dirac-type operator $\D^\pm$, to check that
\begin{eqnarray*}
\D^\mp\Psi^\pm_{a}=\frac{n}{2}H_0\Psi_{a}^\pm
\end{eqnarray*}
which, by Section \ref{Sp-Ch}, only depend on the Riemannian metric and the spin structure on $\Sigma$.

The previous results could be stated as : 
\begin{proposition}\label{IKS}
For any $a\in\mathbb{C}^N$, the sections of $\E\mathbb{H}$ defined by
$$
\Phi_a^\pm:=\left\lbrace
\begin{array}{ll}
\psi_{a}^\pm & \text{if $n$ is odd}\\\\
\Psi_{a}^\pm & \text{if $n$ is even}
\end{array}
\right.
$$
are imaginary Killing spinors on $\E\mathbb{H}$. Moreover, if $\Sigma$ is an oriented hypersurface in $\mathbb{H}^{n+1}$, then $\Phi_a^\pm$ satisfies
\begin{eqnarray*}
\D^\mp\Phi_a^\pm=\frac{n}{2}H_0\Phi_a^\pm
\end{eqnarray*}
and this equation only depends on the Riemannian and spin structures of $\Sigma$.
\end{proposition}
\vspace{0.2cm}

As we will see in the next section, the proof of Theorem \ref{general} relies essentially on (\ref{holoim}). However, as easily seen, this principle depends strongly on spinor data whereas our energy-momentum vector $\mathbf{E}(\Sigma)$ does not. A trick by Wang (p. $285$-$286$ in \cite{Wa1}), generalized by Kwong (Proposition $2.1$ and $2.2$ in \cite{K}), allows to clarify these aspects. Indeed, since for any imaginary Killing $\Phi_a^\pm\in\Gamma(\E\mathbb{H})$ as in Proposition \ref{IKS}, we have $|\Phi^\pm_a|^2=|\psi^\pm_a|^2$, we easily deduce  
\begin{lemma}\label{X-K}
For every imaginary Killing spinor $\Phi^\pm_a\in\Gamma(\E\mathbb{H})$, there exists a vector field $\zeta_a^\pm\in\mathbb{R}^{n+1,1}$ given by
\begin{eqnarray*}
\zeta_{a}^\pm=\mp i\sum_{j=1}^{n+1}\<\gamma_e(\partial_{x_j})a,a\>\partial_{x_j}-|a|^2\partial_t
\end{eqnarray*}
such that
\begin{eqnarray}\label{zet}
|\Phi_{a}^\pm|^2=-2\<{\bf X},\zeta_{a}^\pm\>_{\mathbb{R}^{n+1,1}}.
\end{eqnarray}
Moreover, for every null vector $\zeta=(\zeta_1,\cdots,\zeta_{n+1},1)\in\mathbb{R}^{n+1,1}$, there exists $a\in\mathbb{C}^N$ with $|a|=1$ such that $\zeta=\zeta_{a}^\pm$. Here ${\bf X}=(x_1,...,x_{n+1},t)$ is the position vector field in the Minkowski spacetime and $\gamma_e$ is Clifford multiplication on the Dirac bundle $\E\mathbb{B}$.
\end{lemma}


\subsection{Non-negativity of the quasi-local mass}


In this section, we prove Theorem \ref{general}. More precisely, we have to show that the energy-momentum vector field ${\bf E}(\Sigma)\in\mathbb{R}^{n+1,1}$ defined by (\ref{E-M-V}) is timelike future directed or zero. For this we first recall a characterization of such vector fields given in Lemma $5.2$ of \cite{WY} for $3$-dimensional manifolds but which is easily seen to be true in any dimension. 
\begin{lemma}
A non-zero vector $v=(v_1,\cdots,v_{n+1},w)$ is timelike future directed if and only if $\<v,\zeta\>< 0$ for all $\zeta=(\zeta_1,\cdots,\zeta_{n+1},1)$ with $\sum_{j=1}^{n+1}\zeta^2_j=1$.
\end{lemma}

From this characterization, we first have to prove that $\<{\bf E}(\Sigma),\zeta\>< 0$ for all null vectors $\zeta=(\zeta_1,\cdots,\zeta_{n+1},1)$, that is
\begin{eqnarray}\label{main-int}
\int_\Sigma\Big(\frac{H_0^2-H^2}{H}\Big)\<X,\zeta\>\,d\Sigma< 0
\end{eqnarray}
unless ${\bf E}(\Sigma)=0$. However, Lemma \ref{X-K} ensures that for any null vector $\zeta$ as above there exits $a\in\mathbb{C}^N$ with $|a|=1$ such that $\zeta=\zeta_a^\pm$. Then from (\ref{zet}), the inequality (\ref{main-int}) is equivalent to
\begin{eqnarray}\label{main-int2}
\int_\Sigma\Big(\frac{H_0^2-H^2}{H}\Big)|\Phi_a^\pm|^2\,d\Sigma> 0
\end{eqnarray}
for all $a\in\mathbb{C}^N$ with $|a|=1$. On the other hand, since we assume that $\Sigma$ admits an isometric and isospin immersion $F$ into the hyperbolic space $\mathbb{H}^{n+1}$, by Proposition \ref{IKS} it follows that every imaginary Killing spinor field of the form $\Phi_a^\pm$ induces a solution of the Dirac-type equation $\D^\mp\Phi^\pm_a=\frac{n}{2}H_0\Phi^\pm_a$ on $\Sigma$ (intrinsically to $\Sigma$). Moreover since we assume that $H$ is positive on $\Sigma$, we can apply (\ref{holoim}) to every $\Phi^\pm_a$, to get (\ref{main-int2}).

However, if equality is achieved, it follows from Theorem \ref{OneIm} that the shape operators of
$\Sigma$ with respect to its embedding in $\Omega$ and its immersion in ${\mathbb{H}}^{n+1}_{-k^2}$ are the same so that ${\bf E}(\Sigma)=0$. This implies that ${\bf E}(\Sigma)$ is timelike future directed or zero. 

Suppose now that ${\bf E}(\Sigma)=0$. In this case, we already know that $\Sigma$ is connected and that the second fundamental form $A$ of $\Sigma$ in $\Omega$ agrees with the one of $\Sigma$ in $\mathbb{H}^{n+1}$ denoted by $A_0$. On the other hand, the hyperbolic space $\mathbb{H}^{n+1}$ admits a maximal number of linearly independent imaginary Killing spinor fields, so that we can repeat the
argument in the proof of Theorem \ref{OneIm} for each one of the restrictions to $\Sigma$ of these spinor fields. In this way we obtain a maximal number of imaginary Killing spinor fields defined on $\Omega$. But, according to \cite{B2} (see also \cite{BFGK}), this forces the manifold $\Omega$ to have constant curvature $-1$. Moreover, since $A=A_0$, we can then glue along $\Sigma$ in $\Omega$ the exterior of $\Sigma$ in the hyperbolic space to obtain a smooth complete Riemannian manifold $M$ with constant negative sectional curvature which is isometric to the hyperbolic space at infinity. We easily conclude that $M$ is isometric to $\mathbb{H}^{n+1}$ and then $\Omega$ is isometric to a compact domain of $\mathbb{H}^{n+1}$. Finally we may apply the fundamental theorem of the local theory of hypersurfaces (see Theorem $2.1$ in \cite{abe-koike-yamaguchi}) to deduce that the embedding of $\Sigma$ in $\Omega$ and its immersion in ${\mathbb{H}}^{n+1}$ are congruent. The converse of the equality case in Theorem \ref{general} is straightforward.


\subsection{The two dimensional case}\label{AHset}


In this section, we consider the case $n=2$ which is the most relevant from a physical point of view. More precisely, we propose to define a new notion of local energy-momentum vector by setting
\begin{eqnarray*}
{\bf E}(\Sigma)=\int_{\Sigma}\frac{H^2_0-H^2}{H}{\bf X}\,d\Sigma,
\end{eqnarray*} 
with $(\Sigma,g)$ a topological $2$-sphere, whose Gauss curvature $K>-k^2$ and mean curvature $H>0$, considered as the boundary of a $3$-dimensional compact Riemannian domain $\Omega$ with scalar curvature $R\geq -6k^2$.  Here $H_0$ is the mean curvature of the embedding of $(\Sigma,g)$ into the standard hyperbolic space $\mathbb{H}^{3}_{-k^2}$ (whose existence and uniqueness are proved in \cite{P} and \cite{DCW}). Then it follows easily from Theorem \ref{general}, that if $\Omega$ is not isometric to a domain of $\mathbb{H}^3_{-k^2}$, then the energy-momentum vector ${\bf E}(\Sigma)$ is a timelike future directed vector in $\mathbb{R}^{3,1}$. Moreover, it is zero if and only if $\Omega$ is a domain in the hyperbolic space. For a more precise statement of this result, we refer to Theorem \ref{3D}. \\

We conclude that ${\bf E}(\Sigma)$ has the non negativity and rigidity properties which are needed to define an appropriate notion of quasilocal mass. Another important feature is also required: the limit of ${\bf E}(\Sigma)$ should recover the total energy in the asymptotically hyperbolic case. So let us first recall this setting as well as a notion of total energy for such manifolds defined by Wang \cite{Wa1}. A more general setting is described in \cite{CH}. A complete non compact Riemannian manifold $(M^3,g)$ is asymptotically hyperbolic (AH) if $M$ is the interior of a compact manifold $\overline{M}$ with boundary $\partial\overline{M}$ such that
\begin{enumerate}
\item there is a smooth function $r$ on $\overline{M}$, with $r>0$ on $M$ and $r=0$ on $\partial\overline{M}$, such that $\overline{g}=r^2g$ extends as a smooth Riemannian metric on $\overline{M}$;

\item $|dr|_{\overline{g}}=1$ on $\partial\overline{M}$;

\item $\partial\overline{M}$ is the standard unit sphere $\mathbb{S}^2$;

\item on a collar neighborhood of $\partial\overline{M}$, we have:
\begin{eqnarray*}
g=\sinh^{-2}(r)(dr^2+g_r),
\end{eqnarray*}
where $g_r$ is an $r$-dependent family of metrics on $\mathbb{S}^{2}$ such that:
\begin{eqnarray*}
g_r=g_0+\frac{r^3}{3}h+e.
\end{eqnarray*}
Here $g_0$ is the standard metric on the sphere, $h$ is a smooth symmetric $2$-tensor on $\mathbb{S}^2$ and $e$ is of order $O(r^4)$.
\end{enumerate} 

The following positive mass theorem was proved by Wang:
\begin{theorem}\cite{Wa1}
If $(M^3,g)$ is an (AH) Riemannian manifold such that its scalar curvature satisfies $R\geq -6$, then the energy-momentum vector
\begin{eqnarray*}
\Upsilon=\Big(\int_{\mathbb{S}^2}{\rm tr}_{g_0}(h)d\mathbb{S},\int_{\mathbb{S}^2}{\rm tr}_{g_0}(h)xd\mathbb{S}\Big)\in\mathbb{R}^{3,1}
\end{eqnarray*}
is timelike future directed or zero. It is zero if and only if $(M^3,g)$ is isometric to the hyperbolic space $\mathbb{H}^3$. Here $d\mathbb{S}$ denotes the standard Riemannian measure on the round sphere. 
\end{theorem}

Using a recent work of Kwong and Tam \cite{KT}, we show that our local energy-momentum vector (under some additional technical assumptions) converges to 
the energy-momentum vector $\Upsilon$. In fact, as in \cite{KT}, we assume that the following hold:\\

\noindent {\rm {\bf (A)}\hspace{1.5cm} $\nabla_{\mathbb{S}^2}e$, $\nabla_{\mathbb{S}^2}^2e$, $\nabla_{\mathbb{S}^2}^3e$, $\nabla_{\mathbb{S}^2}^4e$ and $\frac{\partial e}{\partial r}$ are of order $O(r^3)$}\\

\noindent where $\nabla_{\mathbb{S}^2}^k$ is the Levi-Civita connection of order $k$ on tensor fields. Then consider a geodesic sphere $S_r\subset(M,g)$ for $r$ small and let $H$ be its mean curvature. We identify $S_r$ as the standard sphere $\mathbb{S}^2$ with metric $\gamma_r$ induced from $g$. For $r$ small enough, the Gauss curvature of $(S_r,\gamma_r)$ is positive, hence $(S_r,\gamma_r)$ can be isometrically embedded into $\mathbb{H}^3$ by Pogorelov's Theorem. If $X^{(r)}$ denotes this embedding and if $o_r$ is the center of the largest geodesic sphere contained in the interior of $X^{(r)}(S_r)$, then Kwong and Tam prove that we can choose the center of the geodesic balls at a fixed point $o\in\mathbb{H}^3$. In addition to this, they construct isometries $\iota_r$ of $\mathbb{H}^3$ fixing $o$ such that, when $X^{(r)}$ is seen as an embedding of $(S_r,\gamma_r)$ into $\mathbb{R}^{3,1}$ (via $\mathbb{H}^3$), the following expansions hold
$$
\left\lbrace
\begin{array}{rll}
H & = & \cosh r-\frac{1}{4}r^3{\rm tr}_{g_0}(h)+o(r^4)=1+\frac{r^2}{2}-\frac{1}{4}r^3{\rm tr}_{g_0}(h)+o(r^4)\\ \\
H_0 & = & \cosh r+o(r^4)=1+\frac{r^2}{2}+o(r^3)\\ \\
dS_r & = & \big(\frac{1}{\sinh^2 r}+o(r^2)\big)\,d\mathbb{S}=\big(\frac{1}{r^2}+o(\frac{1}{r})\big)\,d\mathbb{S}\\ \\
\iota_r\circ\ X^{(r)}(x) & = & \big(\frac{1}{r}+o(1),\frac{x}{r}+o(\frac{1}{r})\big) 
\end{array}
\right.
$$  
Form these estimates, straightforward calculations show that
\begin{eqnarray*}
\frac{H^2_0-H^2}{H}=\frac{1}{2}{\rm tr}_{g_0}(h)r^3+o(r^3)
\end{eqnarray*}
and 
\begin{eqnarray*}
\iota_r\circ X^{(r)}(x)\,dS_r=\Big(\frac{1}{r^3}+o(\frac{1}{r^2}),\frac{x}{r^3}+o(\frac{1}{r^3})\Big)\,d\mathbb{S},
\end{eqnarray*}
hence,
\begin{theorem}\label{asymptotic}
Let $(M^3,g)$ be a $3$-dimensional (AH) hyperbolic manifold satisfying the assumptions {\bf (A)}, then:
\begin{eqnarray*}
\lim_{r\rightarrow 0}{\bf E}(S_r)=\frac{1}{2}\Upsilon
\end{eqnarray*}
where 
\begin{eqnarray*}
{\bf E}(S_r)=\int_{S_r}\big(\frac{H_0^2-H^2}{H}\big)\iota_r\circ X^{(r)}dS_r.
\end{eqnarray*}
\end{theorem}

As we have seen, the proof of Theorem \ref{3D} makes no use of the Postive Mass Theorem for (AH) manifolds unlike the results of Shi-Tam and Kwong. In fact, combining Theorems \ref{3D} and \ref{asymptotic}, we get an alternative proof of the Positive Mass Theorem of Wang under the additional assumptions ${\rm {\bf (A)}}$.
  


\end{document}